\documentclass[aps,pre,reprint,groupedaddress,showpacs]{revtex4-1}

\usepackage[dvipdfmx]{graphicx}
\usepackage{latexsym}
\usepackage[fleqn]{amsmath}
\usepackage{amssymb}
\usepackage{color}
\usepackage{multirow}

\usepackage{bm}

\begin{document}

\title{Super Generalized Central Limit Theorem\\
--Limit distributions for sums of non-identical random variables
with power-laws--
}


\author{Masaru Shintani}
\email[shintani.masaru.28a@kyoto-u.jp]{}
\author{Ken Umeno}%
\email[umeno.ken.8z@kyoto-u.ac.jp]{}
\affiliation{Department of Applied Mathematics and Physics, Graduate
School of Informatics, Kyoto University, Yoshida Honmachi Sakyo-ku, Kyoto, 606--8501}


\date{Aug 22, 2017}

 \begin{abstract}
 
In nature or societies, the power-law is present ubiquitously,  
  and then it is important to investigate the characteristics of
  power-laws in the recent era of big data.  
In this paper we prove the superposition of non-identical stochastic processes with power-laws
  converges in density to a unique stable distribution.
This property can be used to explain the universality of stable laws such that the sums of
  the logarithmic return of non-identical stock price fluctuations follow stable distributions.
 \end{abstract}

\maketitle


{\it Introduction---.}There are a lot of data that follow the power-laws in the world.
Examples of recent studies include, but are not limited to
the financial market \cite{mandelbrot1997variation,mantegna1994stochastic,mantegna1995scaling,gopikrishnan1998inverse,gabaix2006institutional,denys2016universality,tanaka2016statistical},
the distribution of people's assets \cite{druagulescu2001exponential},
the distribution of waiting times between earthquakes occurring \cite{bak2002unified}
and the dependence of the number of wars on its intensity \cite{roberts1998fractality}.
   It is then important to investigate the general characteristics of
   power-laws. 

In particular, as for the data in the financial market,
Mandelbrot \cite{mandelbrot1997variation} firstly argued that the
distribution of the price fluctuations of
cotton follows a stable law.
Since the 1990's, there has been a controversy as to whether 
the central limit theorem or the generalized central limit theorem
(GCLT) \cite{Kolmogorov} as sums of power-law
distributions can be applied to the data of the logarithmic return of
stock price fluctuations.
In particular, Mantegna and Stanley
argued that the
logarithmic return follows a stable distribution with the power-law
index $\alpha<2$ \cite{mantegna1994stochastic,mantegna1995scaling},
and later they denied their own
 argument by introducing the cubic laws ($\alpha=3$) \cite{gopikrishnan1998inverse}.
 Even recently, some researchers \cite{gabaix2006institutional,denys2016universality,tanaka2016statistical} have argued
 whether a distribution of the logarithmic returns follows power-laws
 with $\alpha>2$ or stable laws with $\alpha<2$.
On the other hand, it is necessary to prepare
 very large data sets to elucidate true tail
behavior of distributions \cite{weron2001levy}. 
 In this respect, the recent study \cite{tanaka2016statistical} showed
 that the large and high-frequency arrowhead data of the Tokyo stock exchange (TSE)
 support stable laws with $1<\alpha<2$.

In this study, we show that the sums of the logarithmic
return of multiple stock price fluctuations follows stable laws, and it
 can be described from a theoretical background.
We will extend
the GCLT to sums of independent non-identical stochastic processes.
We call this Super Generalized Central Limit Theorem (SGCLT). 

{\it Summary of stable distributions and the GCLT---.}A probability density function $S(x;\alpha,\beta,\gamma,\mu)$ of random variable $X$ following a
stable distribution \cite{nolan2003stable} is defined with its characteristic function $\phi(t)$ as:
\begin{eqnarray*}
S(x;\alpha,\beta,\gamma,\mu)=\frac{1}{2\pi}\int_{-\infty}^{\infty}\phi(t)e^{-ixt}{\rm
 d}x,
\end{eqnarray*}
where $\phi(t;\alpha,\beta,\gamma,\mu)$ is expressed as:
\begin{eqnarray*}
&\phi(t)=\exp{\{i\mu
 t-\gamma^\alpha|t|^{\alpha}(1-i\beta {\rm sgn}(t)w(\alpha,t))\}}\\
& w(\alpha,t) =\left\{\begin{aligned}
	       &\tan\left(\pi\alpha/2\right) & \text{if} &\ \  \alpha\neq 1\\
	       &-2/\pi\log|t| & \text{if} &  \ \   \alpha = 1.
	      \end{aligned}\right. 
\end{eqnarray*}
The parameters $\alpha,\beta,\gamma$ and $\mu$ are real constants
 satisfying $0<\alpha\le 2$, $-1\le \beta\le 1$, $\gamma>0$, and denote the
 indices for power-law in stable distributions,
 the skewness, the scale parameter and the location, respectively.
When $\alpha=2$ and $\beta=0$, the probability density
 function obeys a normal distribution.
Note that explicit forms of stable distributions are not known for
 general parameters $\alpha$ and $\beta$ except for a few cases such as
 the Cauchy distribution ($\alpha=1,\beta=0$).

 A stable random variable satisfies the following property for the scale and
the location parameters.
 A random variable $X$ follows $S(\alpha,\beta,\gamma,\mu)$, when
 \begin{eqnarray}\label{eq:scaleandshift}
  X\overset{d}{=}
   \left\{
\begin{aligned}
 &\gamma X_0+\mu &\ \text{if}  & \ \  \alpha\neq 1\\
 &\gamma X_0+\mu+\frac{2}{\pi}\beta\gamma\ln \gamma &\  \text{if}& \ \  \alpha=1,
\end{aligned}
\right.
\end{eqnarray}
where $X_0=S(\alpha,\beta,1,0)$. 
When the random variables $X_j$ satisfy $X_j \sim
S(x;\alpha,\beta_j,\gamma_j,0)$,
 the superposition $Z_n=(X_1+\cdots+ X_n)/n^{\frac{1}{\alpha}}$ of independent random
 variables $\{X_j\}_{j=1,\cdots,n}$ that have {\it different} parameters
 except for $\alpha$
 is also in the stable distribution family as:
\begin{eqnarray}\label{eq:superposition}
 Z_n \sim S(\alpha,\hat{\beta},\hat{\gamma},\hat{\mu}),
\end{eqnarray}
where the parameters $\hat{\beta},\hat{\gamma}$ and $\hat{\mu}$ are expressed as:
\begin{eqnarray*}
  & \displaystyle
  \hat{\beta}=\frac{\sum_{j=1}^{n}\beta_j\gamma_j^\alpha}{\sum_{j=1}^{n}\gamma_j^\alpha},\hat{\gamma}=\left\{\frac{\sum_{j=1}^{n}\gamma_j^\alpha}{n}\right\}^{\frac{1}{\alpha}}\
  \text{and}\\
 & \hat{\mu}=
  \left\{
\begin{array}{lll}
 0 & \text{if} & \alpha\neq 1\\
 -\frac{2\ln n}{n\pi}\sum_{j=1}^n\beta_j\gamma_j & \text{if} & \alpha=1.
\end{array}\right.
\end{eqnarray*}
We can prove this immediately by the use of the characteristic function for
the sums of random variables expressed as the product of their
characteristic functions:
\begin{eqnarray*}
\phi (t;\alpha,\hat{\beta},\hat{\gamma},\hat{\mu})=\prod_{j=1}^n\phi
 \left(t/n^{\frac{1}{\alpha}};\alpha,\beta_j,\gamma_j,0 \right).
\end{eqnarray*}

We focus on the GCLT. 
Let $f$ of $x$ be a probability density function of a random variable
$X$ for $0<\alpha<2$:
\begin{eqnarray}\label{eq:GCLTcondition}
 f(x) \simeq\left\{
  \begin{aligned}
 & c_{+}x^{-(\alpha+1)}  &\text{for} &\ \ x\to \infty\\
 & c_{-}|x|^{-(\alpha+1)} &\text{for} &\ \  x\to -\infty,
   \end{aligned}\right.
\end{eqnarray}
with $c_+,c_- >0$ being real constants.
   Then, according to the GCLT \cite{Kolmogorov}, the
    superposition of independent, identically distributed random
    variables $X_1,\cdots,X_n$ converges in density to a unique stable
   distribution $S(x;\alpha,\beta,\gamma,0)$ for $n \to \infty$, that is
   \begin{eqnarray}\label{eq:GCLT}
    \begin{split}
Y_n&=\frac{\sum_{i=1}^{n}X_i-A_n}{n^{\frac{1}{\alpha}}} \xrightarrow{d}
 S(\alpha,\beta,\gamma,0) \ \text{for} \ n\to \infty,\\
  A_n &= \left\{
\begin{aligned}
&  0&   \text{if} & \ \  0<\alpha<1\\
& n^2  \Im \ln (\varphi_X\left(1/n\right))& 
 \text{if} &\ \  \alpha=1\\
& n\mathbb{E}[X]  &  \text{if} & \ \   1<\alpha<2,
     \end{aligned}\right.
     \end{split}
\end{eqnarray}
 where $\varphi_X$ is a characteristic function of $X$ as the expected value of
  $\exp(itX)$, $\mathbb{E}[X]$ is the expectation value of $X$, $\Im$ is an imaginary
  part of the argument, 
  and parameters $\beta$ and $\gamma$ are expressed as:
\begin{eqnarray*}
\beta&=&\frac{c_{+}-c_{-}}{c_{+}+c_{-}}, \ 
 \gamma=\left\{\frac{\pi(c_{+}+c_{-})}{2\alpha\sin(\frac{\pi\alpha}{2})\Gamma(\alpha)}\right\}^{\frac{1}{\alpha}},
\end{eqnarray*}
with $\Gamma$ being the Gamma function.
When $\alpha= 2$, we obtain $\mu=\int xf(x){\rm d}x$, $\sigma^2=\int
x^2 f(x){\rm d}x$ and the superposition $Y_n$ of the independent, identically distributed
random variables converges in density to a normal distribution:
\begin{eqnarray*}
Y_n&=&\frac{\sum_{i=1}^{n}X_i-n\mu}{\sqrt{n}\sigma} \xrightarrow{d}
 \mathcal{N}(0,1), \ \text{for} \  n\to \infty.
\end{eqnarray*}

{\it Our generalization---.}We consider an extension of this existing theorem for sums of non-identical random variables.
In what follows we assume that the random variables
$\{X_i\}_{i=1,\cdots,n}$ satisfy the following two conditions.

(Condition 1): The random variables $C_{+}>0$, $C_{-}>0$ obey respectively the
 distributions ${\rm P}_{c_+}(c)$, ${\rm P}_{c_-}(c)$, and satisfy
 $\mathbb{E}[C_+]<\infty$, $\mathbb{E}[C_-]<\infty$.

(Condition 2): The probability distribution function $f_i(x)$ of the random variables
 $X_i$ satisfies in $0<\alpha<2$:
 \begin{eqnarray}\label{eq:SGCLTcondition}
   f_i(x) \simeq\left\{
  \begin{aligned}
 & c_{+i}x^{-(\alpha+1)}  &\text{for} &\ \ x\to \infty\\
 & c_{-i}|x|^{-(\alpha+1)} &\text{for} &\ \  x\to -\infty,
   \end{aligned}\right.
\end{eqnarray}
 where $c_{+i}$ and $c_{-i}$ are samples obtained by $C_{+}$ and $C_{-}$.
We emphasize that the probability distribution function may not be obtained even when we integrate $f_i(x)$ over $c_{+i}$ and $c_{-i}$.


 The main claim of this paper is the following generalization of GCLT:
The following superposition $S_n$ of {\it non-identical} random variables with power-laws converges in
density to a {\it unique stable}
   distribution $S(x;\alpha,\beta^*,\gamma^*,0)$ for $n \to \infty$, where
   \begin{eqnarray}\label{eq:SGCLT}
    \begin{aligned}
S_n&=\frac{\sum_{i=1}^{n}X_i-A_n}{n^{\frac{1}{\alpha}}} \xrightarrow{d}
 S(x;\alpha,\beta^*,\gamma^*,0) \ \   \text{for} \  n\to \infty,\\
          A_n &= \left\{
\begin{array}{lll}
  0&   \text{if} & \   0<\alpha<1\\
 n\sum_{i=1}^n  \Im \ln (\varphi_i\left(1/n\right))& 
 \text{if} &\  \alpha=1\\
 \sum_{i=1}^n \mathbb{E}[X_i]&  \text{if} & \   1<\alpha<2,
     \end{array}\right.
     \end{aligned}
\end{eqnarray}
 with $\varphi_i$ being a characteristic function of $X_i$ as the expected value of
  $\exp(itX_i)$,
  and parameters
  $\beta^*,\gamma^*,\beta_i,\gamma_i$ are expressed as:
\begin{eqnarray*}
\beta^*&=&\frac{{\rm E}_{C_+,C_-}[\beta_i \gamma_i^\alpha]}{{\rm
 E}_{C_+,C_-}[\gamma_i^\alpha]}, \ \ \gamma^*=\left\{{\rm E}_{C_+,C_-}[\gamma_i^\alpha]\right\}^{\frac{1}{\alpha}},\\
\beta_i&=&\frac{c_{+i}-c_{-i}}{c_{+i}+c_{-i}}, \ \ 
 \gamma_i=\left\{\frac{\pi(c_{+i}+c_{-i})}{2\alpha\sin(\frac{\pi\alpha}{2})\Gamma(\alpha)}\right\}^{\frac{1}{\alpha}}.
\end{eqnarray*}
Here ${\rm E}_{C_+,C_-}[X]$ denotes the expectation value of $X$ with respect to
random parameter distributions ${\rm P}_{c_{+}}$ and ${\rm P}_{c_{-}}$.

{\it Proof---.}Although the following is not mathematically rigorous,
we give the following intuitive proof. 

The probability distribution function of random variables $\{X_j\}_{j=1,\cdots,N}$
satisfying the Conditions 1-2 is expressed as:
\begin{eqnarray*}
 f_{j}(x)\simeq
  \left\{
   \begin{aligned}
    &c_{+j}x^{-(\alpha+1)} &\text{for} &\  x\to +\infty\\
    &c_{-j}|x|^{-(\alpha+1)} &\text{for} &\  x\to -\infty,
\end{aligned}\right.
\end{eqnarray*}
where $c_{+j}>0$ and $c_{-j}>0$ satisfy
$\mathbb{E}[C_+]>0$ and $\mathbb{E}[C_-]>0$.
The superposition $S_N$ is then defined as:
   \begin{eqnarray*}
    \begin{split}
     S_N&=\frac{\sum_{j=1}^{N}X_j-A_N}{N^{\frac{1}{\alpha}}},\\
       A_N &= \left\{
\begin{array}{lll}
  0&   \text{if} & \  0<\alpha<1\\
 N \sum_{j=1}^N  \Im \ln (\varphi_j\left(1/N\right)) & 
 \text{if} & \  \alpha=1\\
 \sum_{j=1}^N\mathbb{E}[X_j] &  \text{if} & \    1<\alpha<2,
     \end{array}\right.
     \end{split}
\end{eqnarray*}
 where $\varphi_j$ is a characteristic function of $X_j$.
On the other hand, let $N'$ be $M\times N$ with some $M$, and
$\{X_{ij}\}_{i=1,\cdots,M,j=1,\cdots,N}$ be samples given by the same parent to $X_j$ for each $j$.
Then $\{X_{ij}\}_{i=1,\cdots,M,j=1,\cdots,N}$ are independent, identically
distributed for $i=1,\cdots,M$ at a fixed index $j$.
Then, we define the superposition $S_{N'}$ as follows:
\begin{eqnarray*}
\begin{aligned}
  &S_{N'} = \frac{\sum_{i=1}^M\sum_{j=1}^{N}
     X_{ij}-A_{N'}}{N'^{\frac{1}{\alpha}}},\\
    A_{N'} &= \left\{
  \begin{array}{lll}
  0&   \text{if} &\    0<\alpha<1\\
 M^2N \sum_{j=1}^N \left(  \Im \ln (\varphi_{j}\left(1/(MN))\right)\right)  & 
 \text{if} &\   \alpha=1\\
 M\sum_{j=1}^N\mathbb{E}[X_j] &
 \text{if} &\   1<\alpha<2.
 \end{array}\right.
 \end{aligned}
\end{eqnarray*}


Here, we do not consider the convergence of $S_N$ in density for $N \to
\infty$, but consider the superposition $S_{N'}$ for $N'\to
\infty$,
since the superposition $S_N$ will converge to the same limiting distribution of $S_{N'}$
if $S_N$ converges in density.

We focus on the convergence in density of $S_{N'}$ for
$M\to\infty$ and $N\to\infty$ as follows.
About the previous $A_{N'}$ in $S_{N'}$, we express it as
$A_{N'}=\sum_{j=1}^{N}A_{M_{j}}$ with the following $A_{M_j}(j=1,\cdots,N)$,
\begin{eqnarray*}
  A_{M_{j}} &=& \left\{
\begin{aligned}
&  0&   \text{if} & \ \  0<\alpha<1\\
& M^2N \Im \ln (\varphi_{j}\left(1/MN\right))
& \text{if} &  \ \ \alpha=1\\
& M \mathbb{E}[X_{j}] &  \text{if} & \ \  1<\alpha<2.
  \end{aligned}\right.
\end{eqnarray*}
Here, the superposition $S_{N'}$ is described as:
\begin{eqnarray*}
    S_{N'} &=& \frac{\sum_{i=1}^M\sum_{j=1}^{N}
     X_{ij}-A_{N'}}{N'^{\frac{1}{\alpha}}}\\
 &=&   \frac{\frac{\sum_{i=1}^MX_{i1}-A_{M_{1}}}{M^{\frac{1}{\alpha}}}+\cdots+\frac{\sum_{i=1}^MX_{iN}-A_{M_{N}}}{M^{\frac{1}{\alpha}}}}{N^{\frac{1}{\alpha}}}.
\end{eqnarray*}
When $\alpha\neq 1$,
let $Y_{M_j}$ be
the superposition $\left(\sum_{i=1}^MX_{ij}-A_{M_{j}}\right)/{M^{\frac{1}{\alpha}}}$.
Then, $Y_{M_j}$ converges in density to $S(\alpha,\beta_j,\gamma_j,0)$ for $M
\to \infty$ according to the GCLT \eqref{eq:GCLT}, that is
\begin{eqnarray*}
Y_{M_j}=\frac{\sum_{i=1}^MX_{ij}-A_{M_{j}}}{M^{\frac{1}{\alpha}}}
 \overset{d}{\to} S(\alpha,\beta_j,\gamma_j,0) \ \text{for} \ M\to \infty,
\end{eqnarray*}
where $\beta_j$ and $\gamma_j$ are
\begin{eqnarray*}
\beta_j&=&\frac{c_{+j}-c_{-j}}{c_{+j}+c_{-j}}, \ 
 \gamma_j=\left\{\frac{\pi(c_{+j}+c_{-j})}{2\alpha\sin(\frac{\pi\alpha}{2})\Gamma(\alpha)}\right\}^{\frac{1}{\alpha}}.
 \end{eqnarray*}
Thus, with the stable property \eqref{eq:superposition}, we obtain
the convergence of the superposition $S_{N'}$ as follows:
\begin{eqnarray*}
 S_{N'} &=& \frac{\sum_{j=1}^NY_{M_j}}{N^{\frac{1}{\alpha}}}\\
 &\overset{d}{\to}& \frac{\sum_{j=1}^NY_{j}}{N^{\frac{1}{\alpha}}} \
 \ \text{for} \ M\to \infty , (Y_j\sim S(\alpha,\beta_j,\gamma_j,0)) \\
 &\overset{d}{\to}& S(x;\alpha,\beta^*,\gamma^*,0) \ \   \text{for} \ N\to\infty,
\end{eqnarray*}
where $\beta^*$ and $\gamma^*$ are:
\begin{eqnarray*}
\beta^*&=&\lim_{N\to
 \infty}\frac{\sum_{j=1}^{N}\beta_j\gamma_j^\alpha}{\sum_{j=1}^{N}\gamma_j^\alpha}\\
 &=& \lim_{N\to
 \infty}\frac{\frac{1}{N}\sum_{j=1}^{N}\beta_j\gamma_j^\alpha}{\frac{1}{N}\sum_{j=1}^{N}\gamma_j^\alpha}
 = \frac{{\rm E}_{C_+,C_-}[\beta_j\gamma_j^\alpha]}{{\rm E}_{C_+,C_-}[\gamma_j^\alpha]},\\
\gamma^* &=& \lim_{N\to
 \infty}\left\{\frac{\sum_{j=1}^{N}\gamma_j^\alpha}{N}\right\}^{\frac{1}{\alpha}}
 = \left\{{\rm E}_{C_+,C_-}[\gamma_j^\alpha]\right\}^{\frac{1}{\alpha}}.
\end{eqnarray*}
This proves the superposition $S_{N'}$ converges in density to $S(\alpha,\beta^*,\gamma^*,0)$.
Figure \ref{fig:concept} illustrates the concept of this proof.


     \begin{figure}[htb]
      \begin{center}
 \flushleft{ Step(i)}\\
  \begin{tabular}{p{6.2em}|c|p{2em} c p{2em}|c|}\cline{2-2}\cline{6-6}
 & $X_{11}$ & & $\cdots$ & & $X_{1N}$\\ 
{\scriptsize superposition}  & $\vdots$ && $\ddots$ && $\vdots$\\ 
 & $X_{M1}$ && $\cdots$ && $X_{MN}$\\
  \cline{2-2}\cline{6-6}
   \end{tabular}\\[5pt]
\begin{tabular}{p{4.5em}ccccc}
 \ \ \ \ \ {\footnotesize $M\to \infty$} & {\large$\downarrow$}  &&  && {\large$\downarrow$}  \\[4pt]
 & $S(\alpha,\beta_1,\gamma_1,0)$ && $\cdots$ && $S(\alpha,\beta_N,\gamma_N,0)$
 \end{tabular}\\[3pt]
\flushleft{ Step(ii)}   
   \begin{eqnarray*}
\underbrace{
      S(\alpha,\beta_1,\gamma_1,0),\  \cdots
 \     ,S(\alpha,\beta_N,\gamma_N,0)}_{\text{superposition}} \\
   \overset{N\to \infty}{\longrightarrow} S(\alpha,\beta^*,\gamma^*,0)
   \end{eqnarray*}
       \caption{Concept of the convergence (when $\alpha\neq 1$)}
       \label{fig:concept}
   \end{center}
     \end{figure}

As above, the superposition $S_{N'}$ of non-identical stochastic processes
converges in density to a unique stable distribution.
Since the limiting distribution of $S_{N'}$ is the same as that of $S_N$,
$S_N$ also converges to $S(x,\alpha,\beta^*,\gamma^*,0)$.
When $\alpha=1$, this statement does not hold because of dependence between
$M$ and $N$ in $A_{M_j}$, but we find that the limit
distribution of the superposition $S_N$ generally converges in density to
$S(x;\alpha,\beta^*,\gamma^*,0)$ as can be seen in the following numerical examples.


{\it Numerical confirmation---.}As below, we confirm the claim of SGCLT
\eqref{eq:SGCLT} by some numerical experiments.



To verify the main claim numerically,
we use two kinds of test: two-samples Kolmogorov-Smirnov (KS) test \cite{stephens1974edf} and 
two-samples Anderson-Darling (AD) test \cite{anderson1952asymptotic}
with 5\% significance level.
We generate two data by different methods, and see the
$P\mathalpha{-}values$ of both of tests. Then, unless the null hypothesis is rejected,
we judge the two data follow the same distribution.




For the first data, we generate non-identical stochastic processes satisfying
Conditions 1-2, and prepare the superposition obtained in the same way as
\eqref{eq:SGCLT}. For the second data, we generate the random numbers that
follow the stable distribution, where the first data will converge to
the stable distribution according to \eqref{eq:SGCLT}.
Note that we compare the superposition with
 not a cumulative distribution function but
 random numbers obtained from another
 numerical method described below since a cumulative distribution function of a stable
 distribution cannot be expressed explicitly except for a few cases.




For the first data,
let us consider the chaotic dynamical system $x_{n+1}=g(x_n)$, where $g(x)$
is defined
\cite{umeno1998superposition} as follows for $0<\alpha <2$:
\begin{eqnarray*}
 g(x)=\left\{
				   \begin{array}{lll}
				    \frac{1}{\delta_1^2|x|}\left(\frac{|\delta_1
				    x|^{2\alpha}-1}{2}\right)^{1/\alpha}&
				       \text{for} & \ \  x>\frac{1}{\delta_1}\\
				    -\frac{1}{\delta_1\delta_2|x|}\left(\frac{1-|\delta_1
				    x|^{2\alpha}}{2}\right)^{1/\alpha}& 
				     \text{for} & \ \  0<x<\frac{1}{\delta_1}\\
				    \frac{1}{\delta_1\delta_2|x|}\left(\frac{1-|\delta_2
				    x|^{2\alpha}}{2}\right)^{1/\alpha}& 
				     \text{for} & \ \  -\frac{1}{\delta_2}<x<0\\
				    -\frac{1}{\delta_2^2|x|}\left(\frac{|\delta_2
				    x|^{2\alpha}-1}{2}\right)^{1/\alpha}&
				      \text{for}  & \ \  x<-\frac{1}{\delta_2}.
				   \end{array}
				 \right.
\end{eqnarray*}
This mapping has a mixing property
and an ergodic invariant density for almost all initial points $x_0$.
One of the authors (KU) obtained the following explicit asymmetric power-law
distribution as an invariant density \cite{umeno1998superposition}:
\begin{eqnarray*}
  \rho_{\alpha,\delta_1,\delta_2}(x) =\left\{
  \begin{aligned}
 & \frac{\alpha\delta_1^\alpha
   x^{\alpha-1}}{\pi(1+\delta_1^{2\alpha}x^{2\alpha})}  &\text{if} &\ \
   x\ge 0\\
 &  \frac{\alpha\delta_2^\alpha
   |x|^{\alpha-1}}{\pi(1+\delta_2^{2\alpha}|x|^{2\alpha})}  &\text{if}
   &\ \ x<0.
   \end{aligned}\right.
\end{eqnarray*}
This asymmetric distribution behaves as follows for $x \to \pm\infty$:
\begin{eqnarray*}
  \rho_{\alpha,\delta_1,\delta_2}(x) \simeq\left\{
  \begin{aligned}
& \frac{\alpha}{\pi\delta_1^\alpha}x^{-(\alpha+1)} & \text{for} & \ \  x \to +\infty\\
 &  \frac{\alpha}{\pi\delta_2^\alpha}|x|^{-(\alpha+1)} & \text{for} & \ \ x \to -\infty.
      \end{aligned}\right.
\end{eqnarray*}
This is exactly the same expression with the condition of GCLT
\eqref{eq:GCLTcondition} for random variables in $X$.
Then, putting the variables $\delta_1$ and $\delta_2$ be distributed,
we can obtain various {\it different} distributions with the same
power-laws.

We regard the parameters $\delta_{1i}$ and $\delta_{2i}$ as random samples
obtained from $\Delta_1$ and $\Delta_2$, where $\Delta_1$ and $\Delta_2$ obey ${\rm P}_{\delta_1}(\delta)$ and ${\rm P}_{\delta_2}(\delta)$, respectively. These are defined for $\delta>0$ with finite mean.

Then the parameters $c_{+i}$ and $c_{-i}$ are given as
$c_{+i}=\frac{\alpha}{\pi\delta_{1i}^\alpha}$ and
$c_{-i}=\frac{\alpha}{\pi\delta_{2i}^\alpha}$,
and $\mathbb{E}[C_+]<\infty$, $\mathbb{E}[C_-]<\infty$ are also satisfied since $\delta_{1i}, \delta_{2i}$ are not 0 and samples from 
some random variables $\Delta_1$ and $\Delta_2$ with finite mean.
As above, we can get some stochastic processes satisfying the Conditions 1-2.

For the second data,
the random numbers generated with the following procedure
follow a stable distribution \cite{chambers1976method}.
Let $\Theta$ and $\Omega$ be independent random numbers: $\Theta$ uniformly
distributed in $\left(-\frac{\pi}{2},\frac{\pi}{2}\right)$, $\Omega$
exponentially distributed with mean $1$.
In addition, let $R$ be as follows:
\begin{eqnarray*}
R=\left\{
\begin{matrix}
 \frac{\sin\left(\alpha(\theta_0+\Theta)\right)}{(\cos\left(\alpha\theta_0\right)
 \cos\Theta)^{1/\alpha}}\left[\frac{\cos((\alpha-1)\Theta)}{\Omega}\right]^{(1-\alpha)/\alpha}
 & (\alpha\neq 1)\\
\frac{2}{\pi}\left[(\frac{\pi}{2}+\beta\Theta)\tan\Theta-\beta\log\left(\frac{\frac{\pi}{2}\Omega\cos\Theta}{\frac{\pi}{2}+\beta\Theta}\right)\right]
 & (\alpha=1),
\end{matrix} \right.
\end{eqnarray*}
for $0< \alpha \le 2$ where $\theta_0=\arctan(\beta\tan(\pi\alpha/2))$.
Then it follows that $R \sim S(x;\alpha,\beta,1,0)$.
We get arbitrary stable distributions by the use of the property
\eqref{eq:scaleandshift} about the scale parameter and the location.

 \begin{table}[b]
  \begin{center}
\begin{tabular}{|p{1.7em}|p{4em}|p{4em}|p{3.3em}|p{3.5em}||p{4em}|p{4.2em}|} \hline
 \multicolumn{1}{|c|}{$\alpha$} & \multicolumn{1}{|c|}{${\rm P}_{\delta_1}(\delta)$} & \multicolumn{1}{|c|}{${\rm P}_{\delta_2}(\delta)$} & \multicolumn{1}{|c|}{$N$} & \multicolumn{1}{|c||}{$L$} & $P\mathalpha{-}value$ (KS test) & $P\mathalpha{-}value$ (AD test) \\
  \hline \hline
 \multirow{3}*{$0.5$} & $1$(const) & $1$(const) & 10000 & 50000 & 0.122 & 0.074\\ \cline{2-7}
  & ${\rm U}(1,2)$ & ${\rm U}(1,2)$ & 1000 & 100000 &0.561 &0.413 \\ \cline{2-7}
  & ${\rm U}(0.5,1)$ & ${\rm U}(1,2)$ & 1000 & 100000 &0.865
		     &0.546 \\ \hline
 \multirow{3}*{$1$} & $1$ & $1$ & 1000 & 100000 & 0.226 &0.308 \\ \cline{2-7}
  & ${\rm U}(1,2)$ & ${\rm U}(1,2)$ & 1000 & 100000 &0.741 &0.497\\ \cline{2-7}
  & ${\rm U}(0.5,1)$ & ${\rm U}(1,2)$ & 1000 & 100000 &0.659
		     &0.301\\ \hline
 \multirow{3}*{$1.5$} & $1$ & $1$ & 1000 & 100000 & 0.916 & 0.529\\ \cline{2-7}
  & ${\rm U}(1,1.2)$ & ${\rm U}(1,1.2)$ & 10000 & 20000 & 0.768 & 0.548 \\ \cline{2-7}
  & ${\rm U}(0.5,1)$ & ${\rm U}(1.5,2)$ & 10000 & 30000 & 0.108 &
			 0.099 \\ \hline
\end{tabular}
      \caption{$P-values$ of two tests}
   \label{tab:p-value}
   \end{center}
 \end{table}

  \begin{table*}[t]
  \begin{center}
\begin{tabular}{|p{1.7em}|p{1.5em}|p{1.5em}|p{13em}|p{3em}|p{3.5em}|p{4em}|p{4em}|} \hline
 $\alpha$ &\multicolumn{1}{|c|}{${\rm P}_{\delta_1}(\delta)$} & \multicolumn{1}{|c|}{${\rm P}_{\delta_2}(\delta)$} & random variables & N& L& KS test &AD test \\ \hline
\multirow{2}*{$0.5$} &3&1& $ X_i - i/N$ &2000 &10000 & 0.136 &0.110 \\ \cline{2-8}
 &3&1& $ X_i - \  \text{Crand}(0,1)$ &1000 &10000&0.289  & 0.190 \\ \hline
\multirow{2}*{$1$} &3&1& $ X_i - i/N $ &1000& 10000& 0.305 & 0.081 \\ \cline{2-8}
 &3&1& $ X_i - \  \text{Crand}(0,1)$ &2000& 10000& 0.145 & 0.093 \\ \hline
 $1.5$ &3&1& $ X_i - \  \text{Crand}(0,1)$ &1000& 10000& 0.371 & 0.286  \\ \hline
\end{tabular}
      \caption{$P-values$ of two tests}
   \label{tab:p-value2}
   \end{center}
 \end{table*}

 \begin{figure}[b]
   \begin{center}
 	  \includegraphics[width=80mm]{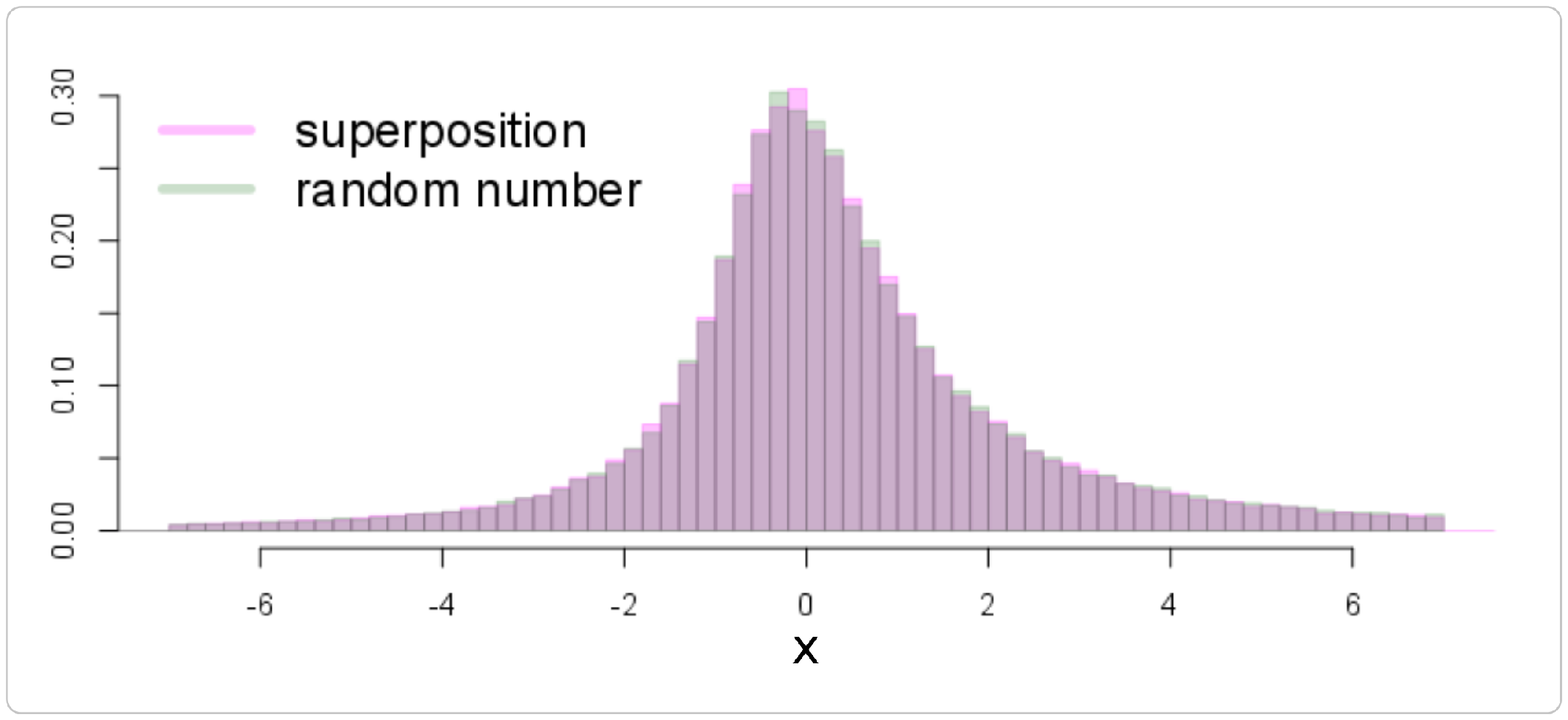}
    \caption{Comparison of two probability densities: the superposition ($N=10^3$, $L=10^5$
    for $\alpha\mathalpha{=}1,\Delta_1\mathalpha{\sim} {\rm U}(0.5,1),\Delta_2\mathalpha{\sim} {\rm U}(1,2)$) and a
    stable distribution ($L\mathalpha=10^5$ for $\alpha\mathalpha{=}1,\beta^*\mathalpha{=}1/3,\gamma^*\mathalpha{=}1$)}
    \label{fig:example}
   \end{center}
 \end{figure}

 With two data obtained accordingly,
 we see whether the superposition
  $S_N=(\sum_{i=1}^NX_i-A_N)/N^{1/\alpha}$
 numerically converges in
 density to a stable distribution $S(x;\alpha,\beta^*,\gamma^*,0)$ or not.
Table \ref{tab:p-value} and \ref{tab:p-value2} show $P\mathalpha{-}values$ of the KS test
and the AD test for each $\alpha,\Delta_1,\Delta_2$.
The constant $L$ is the length of the sequence
and $N$ is the number of sequences used for the superposition.
The meaning of ${\rm U}(a,b)$ is the uniform distribution in $(a,b)$.
Figure \ref{fig:example} illustrates an example of correspondence when $\alpha=1$.
``Crand$(0,1)$'' is the random numbers follow the standard Cauchy distribution.
This case shows that the integral average of the probability distribution function
with the Cauchy distribution is not uniquely determined.


\begin{figure}[h]
    \begin{center}
    \includegraphics[width=8cm]{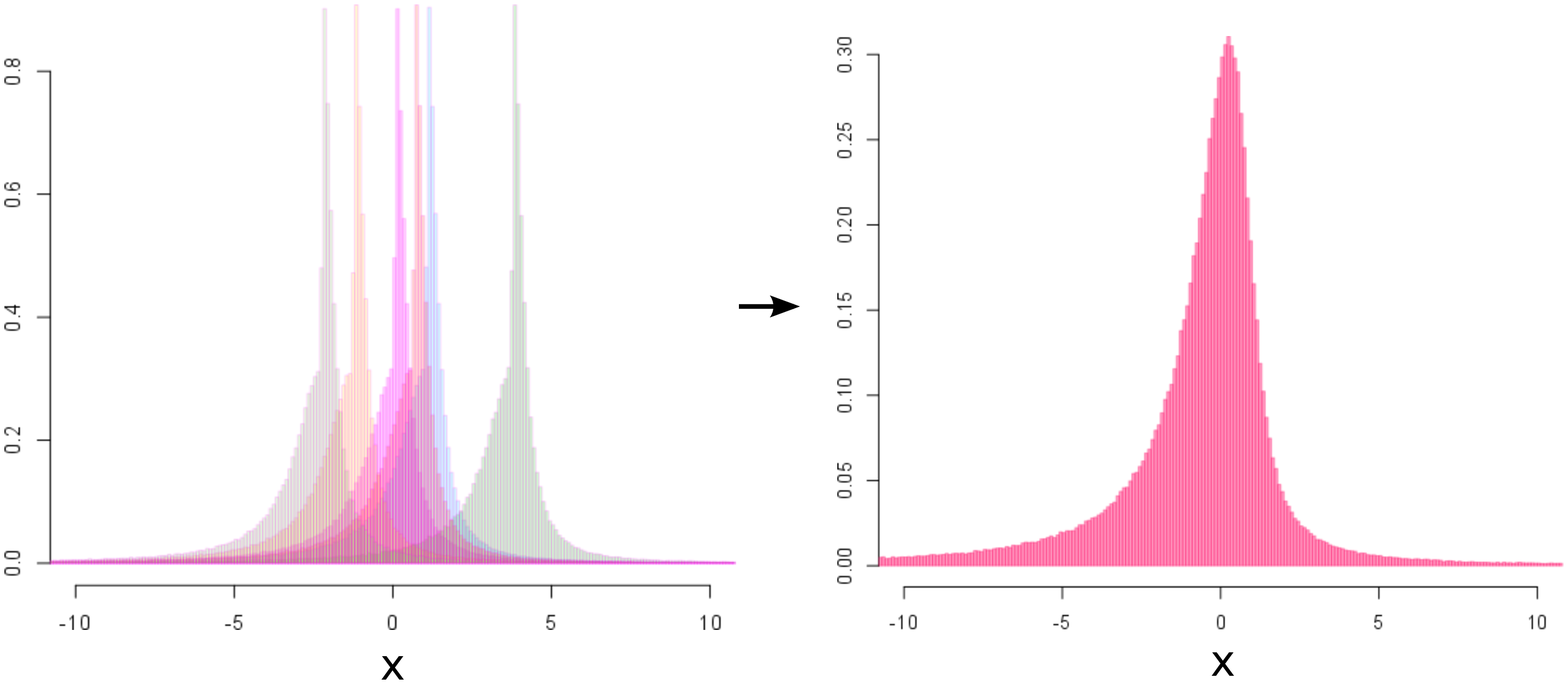}
    \caption{Image of the convergence process: The left figure shows some samples
    of random variables $X_i - \text{Crand}(0,1)$, where $\alpha=1, \delta_1=3, \delta_2=1$. The integration of them does not have an explicit
    expression because of the indefinite mean of the Cauchy
    distribution. However the sum (the right figure) converges to the $S(1,-0.5,2/3,0)$.}
    \label{fig:convergence}
    \end{center}
\end{figure}

As can be seen from Table \ref{tab:p-value} and \ref{tab:p-value2}, we cannot reject the null hypothesis in any
case for $\alpha$.
In other words, the distribution of superposition $S_N$ and the
stable distribution $ S(x; \alpha, \beta^*, \gamma^*,
0)$ are close enough in density according to our SGCLT.

In Figure \ref{fig:convergence}, we can see that
the superposition of non-identical distributed random variables
converges.


{\it Conclusions---.}We have further generalized the GCLT for the sums of independent
{\it non-identical} stochastic processes with the same power-law index
$\alpha$.
Our main claim of SGCLT can have more general applications
since the various type of different power-laws exist in nature.
Thus, our SGCLT can support the argument on the ubiquitous nature of
 stable laws such that the logarithmic return of the
multiple stock price fluctuations follow a stable distribution
 with $1<\alpha<2$ by
 regarding them as the sums of non-identical random variables with power-laws.
Take the data of the stock market as an example.
Then, for the case that the distribution of the logarithmic return of each stock price fluctuation 
have the almost same power-law exponents and different scale parameters $(c_+,c_-)$,
we get some trends or indicators according to this SGCLT.

The authors thank Dr. Shin-itiro Goto (Kyoto University) for stimulating discussions.


\end{document}